\documentclass[12pt]{iopart}
\usepackage{amsfonts}
\usepackage{color}
\usepackage{graphicx}

\newtheorem{theorem}{Theorem}
\newtheorem{lemma}[theorem]{Lemma}
\newtheorem{corollary}[theorem]{Corollary}
\newtheorem{definition}[theorem]{Definition}

\newtheorem{proposition}[theorem]{Proposition}
\newtheorem{remark}[theorem]{Remark}
\def\MM#1{\bi{#1}}
\newcommand{\diff}{\rmd}
\newcommand{\pp}[2]{\frac{\partial #1}{\partial #2}} 
\newcommand{\dede}[2]{\frac{\delta #1}{\delta #2}}
\newcommand{\dd}[2]{\frac{\diff#1}{\diff #2}} 
\newcommand{\eqnref}[1]{(\ref{#1})}
\newcommand{\Vol}{\mathrm{Vol}}
\newcommand{\Id}{\mathrm{Id}}
\begin{document}

\title[The variational particle-mesh method]
{The variational particle-mesh method for matching curves}
\author{C J Cotter}
\address{Department of Aeronautics, Imperial College
  London, SW7 2AZ}
\ead{colin.cotter@imperial.ac.uk}
\begin{abstract}
Diffeomorphic matching (only one of several names for this
technique) is a technique for non-rigid registration of curves
and surfaces in which the curve or surface is embedded in the flow
of a time-series of vector fields. One seeks the flow between two
topologically-equivalent curves or surfaces which minimises some 
metric defined on the vector fields, \emph{i.e.} the flow closest
to the identity in some sense. 

In this paper, we describe a new particle-mesh discretisation for the
evolution of the geodesic flow and the embedded shape. Particle-mesh
algorithms are very natural for this problem because Lagrangian
particles (particles moving with the flow) can represent the movement
of the shape whereas the vector field is Eulerian and hence best
represented on a static mesh. We explain the derivation of the method,
and prove conservation properties: the discrete method has a set of
conserved momenta corresponding to the particle-relabelling symmetry
which converge to conserved quantities in the continuous problem.  We
also introduce a new discretisation for the geometric current matching
condition of (Vaillant and Glaunes, 2005). We illustrate the method
and the derived properties with numerical examples.

\noindent{\it Keywords\/}: symplectic integrators, diffeomorphic shape
matching, geodesics on the diffeomorphism group, EPDiff
\end{abstract}

\ams{65P10,65N21} \pacs{87.57.N-} 
\submitto{\JPA}
\maketitle

\begin{center}
\emph{For Darryl Holm on the occasion of his 60th birthday}
\end{center}

\section{Introduction}
Diffeomorphic matching is a numerical framework for quantifying the
differences between geometric information (such as curves, surfaces,
images or vector fields) using deformations from one geometric object
to another. The geometric object is embedded in the flow generated by
a time-series of vector fields with a chosen norm (such as the
$H^1$-norm) defining the distance along the path between geometric
objects. The computational task is to calculate the velocity fields
which minimise this distance such that one geometric object is mapped
to another. This framework was originally introduced in a series of
papers including \cite{GrMi1998,CaYo01,miller-younes-2001,MiTrYo03}
being extended to match distributions (which can model curves and
surfaces) in \cite{GlTrYo04}, and to match geometric currents (also
for modelling curves and surfaces) in \cite{VaGl2005}. Various
numerical approaches have been proposed for solving the optimisation
problem, either by optimising the functional directly (with an extra
term to penalise flows which do not map close to the target shape) as
described in \cite{BeMiTrYo2003}, or by solving the equations of
motion and shooting for a match between shapes by adjusting the
initial conditions as in \cite{tmt02,McMa2006a,MiMaSh2006}. The main
challenge remains to find a numerical approach which is accurate and
efficient, since the problem of computing the shortest path is a
high-dimensional optimisation problem.

In this paper we introduce a new numerical discretisation for the
diffeomorphic matching problem in the context of matching (although it
can also be used for matching surfaces, images and vector
fields). This method uses a similar approach to the Hamiltonian
Particle-Mesh (HPM) method \cite{FrGoRe2002}, with the difference
being that HPM uses the particle-mesh discretisation to interpolate
density from the particles to the mesh, whereas in this application
the particle-mesh discretisation is used to interpolate momentum
(which takes a central role in the diffeomorphic matching framework).

In section \ref{diffeo} we give a review of the diffeomorphic matching
approach applied to curves in the plane, and establish the notation
which will be used in the other sections. We also discuss the role of
momentum, the conditions used to establish whether the curves have
been matched, and the implications of the particle-relabelling
symmetry satisfied by the equations of motion, as well as the
connection with EPDiff and the Camassa-Holm equation. In section
\ref{particle mesh} we introduce the particle-mesh discretisation and
discuss the discrete symmetries and conservation laws, as well as a
discretisation of the current matching condition and a description of
solution methods. In section \ref{numerics} we illustrate the
properties of the numerical method applied to computing the shortest
path between two test shapes. Finally, in section \ref{summary} we
give a summary and outlook.

\section{Diffeomorphic matching of embedded curves}
\label{diffeo}
In this section we describe the problem of matching one embedded curve
onto another using diffeomorphisms. For simplicity we shall focus on
simple closed curves in the plane although the approach is easily
generalised to other structures.

\subsection{Parameterisations of curves}
We take two simple curves embedded in $\mathbb{R}^2$, $C_A$ and $C_B$,
with parameterisations
\[
\MM{Q}_A(s), \quad \MM{Q}_B(s), \qquad s \in [0,2\pi).
\]
Note that the curves can equally well be represented by
reparameterisations of the curves \emph{i.e.} by
\[
\MM{Q}_A(\eta(s)), \quad \MM{Q}_B(\eta(s)), \qquad s \in [0,2\pi),
\]
where $\eta(s)$ is a diffeomorphism of the circle. Our aim is to find
a matching process which is independent of such reparameterisations.

\subsection{Curves embedded in a flow}
If we take a vector field which defines a fluid flow
$\MM{u}(\MM{x},t)$, and embed a curve in that flow, then the curve
satisfies
\begin{equation}
\label{Q advection}
\pp{}{t}\MM{Q}(t;s) = \MM{u}(\MM{Q}(t;s),t),
\end{equation}
\emph{i.e.} each point of the curve moves with the vector evaluated
at that point.

The aim of the calculation is to search amongst time-series of
vector fields $\MM{u}(\MM{x},t)$, $t\in[0,1]$ such that
\eqnref{Q advection} is satisfied, with the boundary conditions
\[
\MM{Q}(0;s) = \MM{Q}_A(s), \quad \MM{Q}(1;s) = \MM{Q}_B(\eta(s)),
\]
for some (unspecified) reparameterisation $\eta$. If these conditions
are satisfied then we say that $\MM{Q}_A$ is matched onto $\MM{Q}_B$
by the vector field time series $\MM{u}(\MM{x},t)$.

\subsection{Optimisation problem}
We choose a norm for vector fields, such as the $H^n_\alpha$ norm defined by
\[
\|\MM{u}\|^2_{H^n_\alpha}= \int_{\mathbb{R}^2} \MM{u}(\MM{x})\cdot
(1-\alpha^2\nabla^2)^n \MM{u}(\MM{x})\rmd\Vol(\MM{x})
\]
Given a time series $\MM{u}(\MM{x},t)$, we can measure the total
amount of deformation in the flow generated by $\MM{u}$ by the
integral
\begin{equation}
\label{functional}
S(\MM{u}) = \int_0^1\frac{1}{2}\|\MM{u}\|^2\rmd t.
\end{equation}
We wish to find the flow that maps $C_A$ to $C_B$ which is ``nearest
to the identity'' with respect to the choice of norm. This leads to
the following optimal control problem \cite{GlTrYo2006}:
\begin{definition}[Optimisation problem for curve matching]
\label{opti}
Let $\MM{Q}_A,\MM{Q}_B:S^1\to\mathbb{R}^2$ be parameterisations of two
curves $C_A$, $C_B$ in the plane, and let $\|\cdot\|^2$ be a norm for
vector fields in the plane. Then the \emph{distance} between curve $C_A$ 
and $C_B$ is defined to be the minimum over all vector fields $\MM{u}$
of the functional 
\[
\int_0^1\frac{1}{2}\|\MM{u}\|^2 \rmd t
\]
 subject to the following constraints:
\begin{itemize}
\item {\bfseries Dynamic constraint}: $\pp{}{t}\MM{Q}(t;s) =
  \MM{u}(\MM{Q}(t;s),t)$, $\forall t\in[0,1]$.
\item {\bfseries Matching conditions}: $\MM{Q}(0;s) = \MM{Q}_A(s)$,
  $\MM{Q}(1;s) = \MM{Q}_B(\eta(s))$, where $\eta$ is some
  reparameterisation of the circle.
\end{itemize}
\end{definition}

\subsection{Momentum}
\label{momentum}
We can enforce the dynamic constraint by introducing Lagrange
multipliers $\MM{P}(t;s)$, so that the optimisation problem becomes
\begin{eqnarray}
\label{optimisation problem}
\delta \int_0^1 \frac{1}{2}\|\MM{u}\|^2 + \int_{s=0}^{2\pi}
\MM{P}\cdot(\dot{\MM{Q}}-\MM{u}(\MM{Q}))\rmd s \rmd t = 0,
\end{eqnarray}
subject to the matching conditions give above. The Euler-Lagrange
equations then give
\begin{eqnarray}
\label{mom map}
\dede{l}{\MM{u}} & = & \int_{s=0}^{2\pi}\MM{P}\delta(\MM{x}-\MM{Q})\rmd s, \\
\label{Q eqn}
\dot{\MM{Q}} & = & \MM{u}(\MM{Q}), \\
\label{P eqn}
\dot{\MM{P}} & = & -\MM{P}\cdot\nabla\MM{u}(\MM{Q}),
\end{eqnarray}
where $l = \|\MM{u}\|^2/2$. For example, if we choose the $H^n_\alpha$-norm
then 
\[
\dede{l}{\MM{u}} = (1-\alpha^2\nabla^2)^n\MM{u}.
\]
We see that optimal velocity fields take the form
\begin{equation}
\label{opt vels}
\int_{s=0}^{2\pi}\MM{P}G(\MM{x}-\MM{Q})\rmd s, 
\end{equation}
where $G$ is the Green's function associated with the chosen norm,
\emph{e.g.} for the $H^n_\alpha$-norm it is the Green's function for 
the operator $(1-\alpha^2\nabla^2)^n$.

Equation \eqnref{mom map} allows us to write $\MM{u}$ as a function of
$\MM{P}$ and $\MM{Q}$, leading us to notice that equations (\ref{Q
  eqn}-\ref{P eqn}) are canonically Hamiltonian, with Hamiltonian
function given by 
\[
H = l(\MM{u}(\MM{P},\MM{Q})) = \frac{1}{2}\|\MM{u}(\MM{P},\MM{Q})\|^2,
\]
\emph{i.e.}, half the square of the norm of the velocity field written as a
function of $\MM{P}$ and $\MM{Q}$. It is for this reason that we
refer to the Lagrange multiplier $\MM{P}$ as the \emph{momentum}
associated with the curve. This Hamiltonian structure arises from the
fact that the equations have been derived from a variational principle
defined by the extreme points of the functional.  Since the
Hamiltonian for this system is time-independent, it is conserved along
the trajectory, and hence the norm of the velocity is also conserved.

\subsection{Matching condition}
\label{matching condition}
There are a number of different possible ways to pose the matching
condition mathematically. One matching condition which is invariant
under reparameterisations is based on defining a singular vector field
\begin{equation}
\label{singular vector field}
\MM{v}^{\MM{Q}}(\MM{x}) = 
\int_{s=0}^{2\pi}\pp{\MM{Q}}{s}\delta(\MM{x}-\MM{Q})\rmd s.
\end{equation}
We define a functional for these singular vector fields
\begin{eqnarray}
\nonumber
f\left[\MM{v}^{\MM{Q}}\right] &=& \int\MM{v}^{\MM{Q}}
\cdot K * \MM{v}^{\MM{Q}}\rmd\Vol(\MM{x}) \\
&=& \int_{s=0}^{2\pi}\pp{\MM{Q}}{s}\cdot
\int_{s'=0}^{2\pi}\pp{\MM{Q}}{s'}K(\MM{Q}(s)-\MM{Q}(s'))\rmd s\rmd s'.
\label{current matching}
\end{eqnarray}
where $K$ is some smooth kernel function. When $\MM{Q}$ matches
$\MM{Q}_B$ then $f[\MM{v}^{\MM{Q}}-\MM{v}^{\MM{Q}_B}]$ vanishes. This is
called the \emph{current matching} condition \cite{VaGl2005}.


This is a weaker condition that setting the value of $\MM{Q}(s)$ for
each $s$ at time $t=1$, and hence we need to add another boundary
condition to get a unique solution to equations (\ref{mom map}-\ref{P
  eqn}). In \cite{MiTrYo03} it was shown that the solution which
minimises the functional is the one with $\MM{P}$ initially normal to
the curve \emph{i.e.}
\[
\MM{P}\cdot\pp{\MM{Q}}{s} = 0, \qquad \mathrm{at}\quad t=0,
\]
and this extra boundary condition means that equations (\ref{mom
  map}-\ref{P eqn}) have a unique solution.

\subsection{Relabelling}
\label{relabelling}
The optimisation problem in definition \ref{opti} is invariant under
symmetries given by reparameterisations of the curve $\eta$. Noether's
theorem tells us that the equations of motion (\ref{mom map}-\ref{P
  eqn}) have conserved quantities which are generated by these
symmetries. 

To compute the conserved quantities we compute the infinitesimal
generators of these symmetries
\[
\delta\MM{Q}(s;t) = \pp{\MM{Q}}{s}\cdot\xi(s)
\]
where $\xi(s)$ is a vector field on the circle. The cotangent
lift of this infinitesimal symmetry is
\[
\delta\left(\MM{Q}(s;t),\MM{P}(s;t)\right) = 
\left(\pp{\MM{Q}}{s}\cdot\xi(s), -\MM{P}(s;t)\cdot
\pp{\MM{Q}}{s}\cdot\xi(s)\right)
\]
which is a Hamiltonian flow generated by the Hamiltonian functional
\[
h_{\xi} = \int_{s=0}^{2\pi}\MM{P}(s)\cdot\pp{\MM{Q}}{s}\cdot\xi(s)\rmd s.
\]
By Noether's theorem, $h_{\xi}$ is conserved along solutions of
\eqnref{mom map}-\eqnref{P eqn}, and since $\xi(s)$ is arbitrary, we
see that
\[
\MM{P}(s)\cdot\pp{\MM{Q}}{s}
\]
is conserved along solution trajectories, for each $s$. In particular,
this means that if the momentum is initially normal to the curve 
\emph{i.e.} 
\[
\MM{P}(s;0)\cdot\pp{\MM{Q}}{s}(s;0)=0,
\]
then the momentum is normal to the curve for all values of $t$ along
the solution. Therefore, optimal velocity fields take the form of
equation \eqnref{opt vels} with the added constraint that the momentum
$\MM{P}$ is normal to the curve.

\subsection{EPDiff}
\label{epdiff}
As described in \cite{CoHo2007,CoHoHy2007}, because the dynamic
constraint is given by a Lie algebra action of velocity fields on
embedded curves, it is possible to eliminate $\MM{P}$ and $\MM{Q}$ by
taking the time derivative of equation \eqnref{mom map} and making use
of equations (\ref{Q eqn}-\ref{P eqn}). This leads to a PDE defined on
the whole of $\mathbb{R}^2$ given by
\[
\MM{m}_t + \nabla\cdot(\MM{u}\MM{m}) + (\nabla\MM{u})^T\MM{m} = 0,
\qquad \MM{m} = \dede{l}{\MM{u}}.
\]
This is the Euler-Poinc\'are equation on the diffeomorphism group,
abbreviated as EPDiff \cite{HoMa2004,HoMaRa1998}, which is the
equation for geodesic flow on the diffeomorphism group. Although we do
not explicitly solve EPDiff during the computation of the optimal
flow, it is useful to understand the computed solutions as singular
solutions of EPDiff. In one dimension and with the $H^1_\alpha$-norm,
EPDiff becomes the Camassa-Holm equation \cite{CaHo1993} which is
completely integrable with singular soliton solutions.

\section{Particle-mesh discretisation}
\label{particle mesh}
We take a \emph{discrete mechanics and optimal control}
\cite{JuMaOb2005} approach by applying the discretisation directly to
the functional \eqnref{optimisation problem} and deriving the
resulting equations. The discretisation used is a particle-mesh method
with:
\begin{itemize}
\item the vector fields discretised on a fixed, finite \emph{mesh}, and
\item the curve discretised as a set of moving points.
\end{itemize}
The principle benefits of this discretisation are that different norms
can be defined on the mesh without the need to calculate Green's
functions. With large numbers of points the method becomes very
efficient as particle momenta can be interpolated to the mesh, then
the norm operator can be inverted, then the velocity values can be
interpolated back to the particle positions.

\subsection{Mesh discretisation of vector fields}
We take a fixed set of points on a mesh $\{\MM{x}_k\}_{k=1}^{n_g}$
(for the numerical examples in this paper we used an equispaced square
mesh) and give each mesh point $\MM{x}_k$ a vector
$\MM{u}_k$.  We interpolate from the set $\{\MM{u}_k\}_{k=1}^{n_g}$ of
vectors to a general point $\MM{x}$ in the plane using a linear
interpolation
\begin{equation}
\label{u formula}
\MM{u}(\MM{x}) = \sum_{k=1}^{n_g}\MM{u}_k\psi_k(\MM{x}).
\end{equation}
For the examples in this paper we set $\psi_k$ to be a tensor product
of cubic B-spline functions centred on $\MM{x}=\MM{x}_k$.

\begin{remark}
The set of vectors on the mesh generate a finite dimensional subspace
of the infinite dimensional space of vector fields. However, the
finite dimensional subspace is not closed under the Lie bracket for
vector fields which means we will not be able to obtain a discrete EPDiff
equation from the discrete equations of motion by eliminating $\MM{Q}$
and $\MM{P}$ as in section \eqnref{epdiff}.
\end{remark}

Once we have this discrete representation of the vector fields we can
define a discretised norm using standard mesh methods (such as finite
difference, finite volume \emph{etc.}). For the the examples in this
paper we took periodic boundary conditions, and used a spectral
discretisation of the $(1-\alpha^2\nabla^2)^n$ operator using discrete
Fourier transforms, and a simple Riemann sum for the integration
(which gives spectral accuracy in this case). Since the Green's
functions decay exponentially over the lengthscale $\alpha$, the
boundary conditions should not affect the solution as long as the
curve is sufficiently far from the boundary. Other boundary conditions
are possible if, for example, a finite-difference method is used to
discretise the norm operator.

\subsection{Particle discretisation of curves}
We replace the parameterised curve $\MM{Q}(s)$ by a finite set of
points in the plane $\{\MM{Q}_\beta\}_{\beta=1}^{n_p}$. The equation
\eqnref{Q advection} gets replaced by the semi-discrete (continuous
time/discrete space) equation
\begin{equation}
\label{semidiscrete Q advection}
\dot{\MM{Q}}_\beta = \sum_{k=1}^{n_g}\MM{u}_k\psi_k(\MM{Q}_\beta).
\end{equation}
We also have to choose a discretisation of integration around the
loop; again for the examples in this paper we use a Riemann sum.

\subsection{Semi-discrete functional and equations of motion}
After the particle-mesh discretisation we obtain a semi-discrete
functional by integrating the discrete vector field norm from $t=0$ to
$t=1$ and introducing Lagrange multipliers to enforce the constraint
\eqnref{semidiscrete Q advection} to get
\begin{equation}
\label{semidiscrete action}
\int_0^1\frac{1}{2} \|\MM{u}\|^2_g + \sum_\beta\MM{P}_\beta\cdot
\left(\dot{\MM{Q}}_\beta - \sum_{k=1}^{n_g}\MM{u}_k\psi_k(\MM{Q}_\beta)
\right) \rmd t.
\end{equation}
The Euler-Lagrange equations for the optimal values of this functional
are
\begin{eqnarray}
\label{semidiscrete mom map}
\pp{}{\MM{u_k}}\frac{1}{2}\|\MM{u}\|^2_g 
&=& \sum_\beta\MM{P}_\beta\psi_k(\MM{Q}_\beta), \\
\label{semidiscrete Q eqn}
\dot{\MM{Q}}_\beta & = &  \sum_{k=1}^{n_g}\MM{u}_k\psi_k(\MM{Q}_\beta), \\
\label{semidiscrete P eqn}
\dot{\MM{P}}_\beta & = & -\MM{P}_\beta\cdot\sum_k\MM{u}_k\nabla
\psi_k(\MM{Q}_\beta).
\end{eqnarray}
Once again, equation \eqnref{semidiscrete mom map} allows us to write
$\MM{u}$ as a function of $\MM{P}$ and $\MM{Q}$, leading us to notice
that equations (\ref{semidiscrete Q eqn}-\ref{semidiscrete P eqn}) are
canonically Hamiltonian, with Hamiltonian function given by
$\|\MM{u}(\MM{P},\MM{Q})\|^2_g/2$.

\subsection{Time discretisation}
\label{time discretisation}
Equation \eqnref{semidiscrete Q advection} can be discretised using
any one-step method, such as a Runge-Kutta method. For simplicity,
in this paper we use the forward Euler discretisation, but the derivation
of the equations is very similar for other methods. The equation becomes
\[
\MM{Q}^{n+1}_\beta = \MM{Q}^n_\beta + 
\Delta t\sum_k\MM{u}^{n+1}\psi_k(\MM{Q}^n).
\]
The discrete functional becomes
\begin{equation}
\label{discrete functional}
\sum_{n=1}^N\left(
\Delta t\frac{1}{2}\|\MM{u}\|^2_g + \sum_\beta\MM{P}^{n+1}_\beta\cdot
\left(\MM{Q}^{n+1}_\beta -\MM{Q}^n_\beta 
- \Delta t\sum_{k=1}^{n_g}\MM{u}_k^{n+1}\psi_k(\MM{Q}_\beta^n)
\right)\right), 
\end{equation}
where $\Delta t = 1/N$, and the discrete Euler-Lagrange equations are
\begin{eqnarray}
\label{discrete mom map}
\pp{}{\MM{u_k}}\frac{1}{2}\|\MM{u}^{n+1}\|^2_g 
&=& \sum_\beta\MM{P}_\beta^{n+1}\psi_k(\MM{Q}_\beta^n), \\
\label{discrete Q eqn}
\MM{Q}^{n+1}_\beta & = & \MM{Q}^n_\beta
+ \Delta t \sum_k\MM{u}^{n+1}_k\psi_k(\MM{Q}_\beta^n), \\
\label{discrete P eqn}
\MM{P}^{n+1}_\beta & = & \MM{P}^n_\beta
- \Delta t\MM{P}^{n+1}\cdot\sum_k\MM{u}_k^{n+1}\nabla\psi_k(\MM{Q}_\beta^n).
\end{eqnarray}

These equations provide a variational integrator \cite{LeMaOrWe2003}
for the semi-discrete equations (\ref{semidiscrete mom
  map}-\ref{semidiscrete P eqn}). Hence, they give a symplectic
integrator (see \cite{LeRe2005} for a survey of these methods) for the
semi-discrete equations in Hamiltonian form. In this particular case
we obtain the first-order symplectic Euler method; higher-order
partitioned Runge-Kutta methods can be obtained by discretising
\eqnref{semidiscrete Q advection} using any Runge-Kutta method. 

\subsection{Discrete symmetries}
\label{discrete relabelling}
In this section we discuss the properties of applying a variational
integrator to equations \eqnref{semidiscrete mom
  map}-\eqnref{semidiscrete P eqn}, and their possible benefits for
the problem of matching curves.

 The properties that make variational integrators the best choice for
 long integrations are:
\begin{itemize}
\item {\bfseries Modified Hamiltonian}: if the Hamiltonian
  $H(\MM{P},\MM{Q})$ is analytic then it is possible to use
  \emph{backward error analysis} to find a
  modified Hamiltonian
\[
\tilde{H} = H + \Delta t^p \Delta H(\MM{P},\MM{Q};\Delta t)
\]
which is conserved over times
\[
|t|\leq c_0\exp^{c_1/\Delta t}.
\]
See \cite{LeRe2005} for a summary and references.
\item {\bfseries Conserved momenta}: if the continuous-time
  Euler-Lagrange equations have conserved momenta associated with
  symmetries of the Lagrangian, then discrete Euler-Lagrange equations
  also conserve these momenta \cite{LeMaOrWe2003}.
\end{itemize}
\subsection{Modified Hamiltonian}

In the problem of matching curves, the Hamiltonian is half the squared
norm for vector fields, which is conserved along trajectories as
described in section \ref{momentum}. If a symplectic integrator is
used then backward error analysis guarantees that the Hamiltonian will
be approximately conserved for long times (although it should be noted
that no theory exists for the piecewise-cubic functions used in the
examples of this paper). It remains an open question as to whether the
conservation of the Hamiltonian is important for these problems on
finite-time intervals, since they are on short time intervals and so
the value of the Hamiltonian will also be approximately preserved by
variable step-size methods using error control.

The main cost of using variational integrators is that they are
implicit for this type of system where the Hamiltonian is a function
of $\MM{P}$ and $\MM{Q}$, but during the optimisation algorithm we are
able to store values along the trajectory from previous integrations
which can be used as initial guesses for solving the implicit
equations using Newton iteration.  The variational integrator allows
large timesteps to be taken in that case (although of course this is
the case for any implicit method, symplectic or otherwise).
\subsection{Momentum conservation}
As described in section \ref{relabelling}, the equations (\ref{mom
  map}-\ref{P eqn}) have a symmetry under reparameterisation of the
curve which has an associated conserved momentum
$\MM{P}\cdot\partial\MM{Q}/\partial s$. This means that for the
optimal solution, $\MM{P}$ remains normal to the curve along the
trajectory. Since we have discretised the curve we have broken that
symmetry, but as we shall show in this section, it is still possible
to recover a sense in which $\MM{P}$ is normal to the curve along the
solution.

 Equation \eqnref{u formula} defines a vector field over the whole
 space parameterised by the vectors $\{\MM{u}_k\}_{k=1}^{n_g}$
 associated with the mesh points. This means that we can choose a
 parameterised curve $\MM{Q}(s)$ which passes through our discrete
 points $\{\MM{Q}_\beta\}_{\beta=1}^{n_p}$ in sequence and follow its
 evolution along the flow by solving
 \begin{equation}
\label{continuous curve}
 \pp{}{t}\MM{Q}(s) = \sum_k\MM{u}_k\psi_k(\MM{Q}(s)).
 \end{equation}
 More generally, we can follow how the flow generated by the vector
fields evolves over the whole space
\[
\pp{}{t}\MM{g}(\MM{x},t) = \sum_k\MM{u}_k\psi_k(\MM{g}(\MM{x},t)),
\]
where $g(\MM{x},t)$ is the flow map taking points from their position
at time $0$ to their position at time $t$ ($\MM{x}$ plays the role of
a coordinate on ``label space'' here), and we can follow the Jacobian
of this map
\[
\pp{}{t}\pp{\MM{g}(\MM{x},t)}{\MM{x}}
 = \sum_k\MM{u}_k\nabla\psi_k(\MM{g}(\MM{x},t))\cdot
\pp{\MM{g}(\MM{x},t)}{\MM{x}}.
\]
In particular, we can evaluate this equation at each of our
discrete points $\MM{Q}_\beta$:
\begin{equation}
\dot{J}_\beta = \sum_k\MM{u}_k\nabla\psi_k(\MM{Q}_\beta)\cdot
J_\beta, \qquad J_\beta = \pp{\MM{g}}{\MM{x}}(\MM{Q}_\beta(0),t),
\qquad J_\beta(0) = \Id.
\label{J evolution}
\end{equation}

A variation $\delta\MM{Q}_\beta(0)$ in the initial conditions
$\MM{Q}_\beta(0)$ leads to a variation in the entire trajectory given
by
\begin{equation}
\label{delta Q}
\delta\MM{Q}_\beta(t) = J_\beta(t)\delta\MM{Q}_\beta(0).
\end{equation}
This variation generates a symmetry of the equations, as shown in
the following lemma.
\begin{lemma}
  The infinitesimal transformation given by equations \eqnref{J
    evolution}-\eqnref{delta Q} together with
\[
\delta\MM{P}_\beta=0, \quad \beta=1,\ldots,n_p, \qquad
\delta\MM{u}_k=0, \quad k=1,\ldots,n_k,
\]
is a symmetry of equations
  \eqnref{semidiscrete mom map}-\eqnref{semidiscrete P eqn}.
\end{lemma}
\noindent\emph{Proof.} Since the equations have been derived from an
action principle, we
  simply need to show that the infinitesimal transformation causes the
  action \eqnref{semidiscrete action} to vanish. If we apply the
  transformation to the integrand (the Lagrangian) then we obtain
\begin{eqnarray*}
\delta L & = & \delta\left(\frac{1}{2}\|\MM{u}\|^2_g + \sum_\beta\MM{P}_\beta\cdot
\left(\dot{\MM{Q}}_\beta - \sum_{k=1}^{n_g}\MM{u}_k\psi_k(\MM{Q}_\beta)\right)
\right), \\
& = & \sum_\beta\MM{P}_\beta\cdot
\left(\delta\dot{\MM{Q}}_\beta - 
\sum_{k=1}^{n_g}\MM{u}_k\nabla\psi_k(\MM{Q}_\beta)\cdot\delta\MM{Q}_\beta
\right), \\
& = & \sum_\beta\MM{P}_\beta\cdot
\left(\dot{J}_\beta\cdot\delta\MM{Q}_\beta(0) - 
\sum_{k=1}^{n_g}\MM{u}_k\nabla\psi_k(\MM{Q}_\beta)\cdot J_\beta\delta\MM{Q}_\beta
(0)
\right), \\
& = & \sum_\beta\MM{P}_\beta\cdot
\left(\sum_k\MM{u}_k\nabla\psi_k(\MM{Q}_\beta)\cdot
{J}_\beta\cdot\delta\MM{Q}_\beta(0) - 
\sum_{k=1}^{n_g}\MM{u}_k\nabla\psi_k(\MM{Q}_\beta)\cdot J_\beta\delta\MM{Q}_\beta
(0)
\right), \\
& = & 0,
\end{eqnarray*}
and hence the result.

\noindent $\square$

This symmetry of the equations has an associated conserved momenta following
Noether's theorem, as described in the following proposition.
\begin{proposition}
For each $\beta = 1,\ldots, n_p$, the quantity $J^T_\beta\MM{P}_\beta$
is conserved.
\end{proposition}
\noindent\emph{Proof.} Applying the symmetry to the action principle 
and substituting the equations of motion \eqnref{semidiscrete mom
  map}-\eqnref{semidiscrete P eqn} gives
\begin{eqnarray*}
0 &=& \delta
\int_0^1 \frac{1}{2}\|\MM{u}\|^2_g + \sum_\beta\MM{P}_\beta\cdot
\left(\dot{\MM{Q}}_\beta - \sum_{k=1}^{n_g}\MM{u}_k\psi_k(\MM{Q}_\beta)
\right) \rmd t \\
& = & \int_0^1
\sum_\beta\left(
\dot{\MM{Q}}_\beta -  \sum_{k=1}^{n_g}\MM{u}_k\psi_k(\MM{Q}_\beta)
\right)\cdot\delta\MM{Q}_\beta \\
& & \quad + 
\dd{}{t}\left(\sum_\beta\MM{P}_\beta J_\beta\cdot\delta\MM{Q}_\beta(0)\right)
\rmd t \\
&= &\sum_\beta\left(\MM{P}_\beta(t_1)J_\beta(t_1)-\MM{P}_\beta(t_0)
J_\beta(t_0)\right)\cdot\delta\MM{Q}_\beta(0).
\end{eqnarray*}
We obtain the result since $\delta\MM{Q}_\beta(0)$ is an arbitrary
vector.

\noindent $\square$

As described in \cite{LeMaOrWe2003}, these conservation laws satisfied
by the semi-discrete equations will be preserved by numerical methods
provided that they are derived from a discrete variational principle
\emph{i.e.}  following the framework described in section \ref{time
  discretisation}.

If we choose a continuous curve through our set of points that moves
with the flow according to equation \eqnref{continuous curve}, then an
approximation to $\partial\MM{Q}(s_\beta)/\partial s\cdot\rmd s$ is
given at time $t$ by
\[
J_\beta\Delta\MM{Q}_\beta(0), \qquad \beta=1,\ldots,n_p,
\]
where $\Delta\MM{Q}_\beta(0)$ is an approximation to
$\partial\MM{Q}(s_\beta)/\partial s\cdot\rmd s$ at time $t=0$. This leads
to the following corollary:
\begin{corollary}
If $\MM{P}_\beta$ is chosen to be normal to $\Delta\MM{Q}_\beta(0)$ at time
$t=0$ for all $\beta$ \emph{i.e.} $\MM{P}$ is initially normal to the curve,
then $\MM{P}$ is normal to $J_\beta\Delta\MM{Q}_\beta(0)$ for all $t$ along
the flow, \emph{i.e.} $\MM{P}$ stays normal to the curve.
\end{corollary}
\noindent\emph{Proof.} For each $\beta=1,\ldots,n_p$, the component of
$\MM{P}$ tangential to the shape is
\begin{eqnarray*}
\MM{P}_\beta(t)\cdot J_\beta(t)\Delta\MM{Q}_\beta(0)
&=& \MM{P}_\beta(0)\cdot J_\beta(0)\Delta\MM{Q}_\beta(0) \\
&=& \MM{P}_\beta(0)\cdot \Id\Delta\MM{Q}_\beta(0) \\
&=& \MM{P}_\beta(0)\Delta\MM{Q}_\beta(0)=0.
\end{eqnarray*}

\noindent $\square$

If the equations of motion for $\MM{Q}$ and $\MM{P}$ are discretised
using a variational integrator then this property is preserved,
provided that the discrete equation for the evolution of $J$ is
defined as the gradient of the evolution equation for $\MM{Q}$, since
then it generates a symmetry of the discrete variational principle
just as in the continuous time case.

For example, the discrete equation for the gradient of the
time-discrete flow for $\MM{Q}$ obtained from the symplectic Euler
method is
\begin{equation}
\label{J discrete}
J_\beta^{n+1} = J_\beta^n + \Delta t\sum_k\MM{u}_k^{n+1}\nabla\psi(\MM{Q}_\beta^n)
J_\beta^n.
\end{equation}
\begin{proposition}
For each $\beta = 1,\ldots, n_p$, the quantity $J^T_\beta\MM{P}_\beta$
is conserved when the equations are integrated using the symplectic
Euler method, and $J$ is obtained from equation \ref{J discrete}.
\end{proposition}
\noindent\emph{Proof.} 
\begin{eqnarray*}
\delta A & = & \delta
\sum_{n=1}^N\left(
\frac{1}{2}\Delta t\|\MM{u}\|^2_g + \sum_\beta\MM{P}^{n+1}_\beta\cdot
\left(\MM{Q}^{n+1}_\beta -\MM{Q}^n_\beta 
- \Delta t\sum_{k=1}^{n_g}\MM{u}_k^{n+1}\psi_k(\MM{Q}_\beta^n)
\right)
\right), \\
& = & \sum_{n=1}^N\left(\sum_\beta\MM{P}^{n+1}_\beta\cdot
\left(\delta\MM{Q}^{n+1}_\beta -\delta\MM{Q}^n_\beta -
\sum_{k=1}^{n_g}\MM{u}_k^{n+1}\nabla\psi_k(\MM{Q}_\beta^n)\cdot
\delta\MM{Q}_\beta^n
\right)\right), \\
& = & \sum_{n=1}^N\sum_\beta\MM{P}^{n+1}_\beta\cdot
\left((J_\beta^{n+1}-J_\beta^n)\cdot\delta\MM{Q}_\beta^0-
\sum_{k=1}^{n_g}\MM{u}_k^{n+1}
\nabla\psi_k(\MM{Q}_\beta^n)\cdot J_\beta^n\delta\MM{Q}_\beta^0
\right), \\
& = & \sum_\beta\MM{P}_\beta^{n+1}\cdot
\left(\sum_k\MM{u}_k^{n+1}\nabla\psi_k(\MM{Q}_\beta^n)\cdot
{J}_\beta^n\cdot\delta\MM{Q}_\beta^0 - 
\sum_{k=1}^{n_g}\MM{u}_k^{n+1}
\nabla\psi_k(\MM{Q}_\beta^n)\cdot J_\beta^n\delta\MM{Q}_\beta^0
\right), \\
& = & 0,
\end{eqnarray*}
Applying the symmetry to the discrete 
action principle \ref{discrete functional}
and substituting the equations of motion \eqnref{semidiscrete mom
  map}-\eqnref{semidiscrete P eqn} gives
\begin{eqnarray*}
0 &=& \sum_\beta\left(\MM{P}_\beta^{n+1}J_\beta^{n+1}-\MM{P}_\beta^n
J_\beta^n\right)\cdot\delta\MM{Q}_\beta^0.
\end{eqnarray*}
We obtain the result since $\delta\MM{Q}_\beta^0$ is an arbitrary
vector.
\noindent $\square$

We can verify the discrete conservation of $\MM{P}_\beta J_\beta$
directly since
\begin{eqnarray*}
\MM{P}^n_\beta J_\beta^n & = & \MM{P}^{n+1}_\beta
(\Id + \sum_k\MM{u}_k^{n+1}\nabla\psi_k(\MM{Q}^n_\beta))J_\beta^n \\
& = & \MM{P}^{n+1}_\beta J^{n+1}_\beta.
\end{eqnarray*}
Similar results can be obtained when higher-order variational
integrators are used by discretising the $\MM{Q}$ equation using a
Runge-Kutta method.

As discussed in section \ref{matching condition}, the norm-minimising 
solution has momentum normal to the shape along the whole trajectory.
Good preservation of the tangential component of momentum is important
because it means that it is only necessary to constrain the momentum
to be normal to the shape in the initial conditions when a shooting
algorithm is used (discussed in section \ref{solution methods}). 

\subsection{Discrete matching condition}
\label{discrete matching condition}
We can also apply a particle-mesh discretisation to the matching
condition described in section \ref{matching condition}. The
particle-mesh representation of equation \eqnref{singular vector field}
is 
\[
\MM{v}^{\MM{Q}}_k = \sum_\beta\Delta\MM{Q}_\beta\psi_k(\MM{Q}_\beta),
\]
where $\Delta\MM{Q}_\beta$ is an approximation to $\partial\MM{Q}(s)/
\partial s\,\rmd s$ at $s=s_\beta$. For example, a simple finite
difference approximation gives
\[
\MM{v}^{\MM{Q}}_k = 
\sum_\beta(\MM{Q}_\beta-\MM{Q}_{\beta-1})\psi_k(\MM{Q}_\beta).
\]
This is the approximation which we will use in the computed numerical
examples used in this paper. 

Once the singular vector field is evaluated on the mesh we can apply
standard mesh discretisation methods to compute the functional
\eqnref{current matching}:
\[
\hat{f}[\MM{v}^{\MM{Q}}] = \sum_{kl}K_{kl}\MM{v}^{\MM{Q}}_k\MM{v}^{\MM{Q}}_l
\]
where $K_{kl}$ is a discretisation of the kernel operator $K$. For the
computed numerical examples used in this paper, the kernel operator
used was the inverse operator $(1-\alpha^2\nabla^2)^{-2}$ discretised
using discrete Fourier transforms.

This results in a numerical discretisation of the functional used for
the current matching condition. After numerical discretisation the
functional will not have a minimum at zero any more, and we must aim
to minimise this functional rather than find the zero to achieve
matching. 

\subsection{Efficiency}
In \cite{VaGl2005}, a mesh-free method was applied to the current
matching problem. The main cost in this type of method is in summing
up Green's functions on each of the points on the curve to calculate
their velocity, and the value of the matching condition. They employed
a fast multipoles algorithm \cite{GrSt1991} which has a computational
cost of size $\mathcal{O}(N\log N)$ (where $N$ is the total number of
particles), although the multiplicative constant in the scaling can be
relatively big, requiring larger $N$ before benefits are seen.

For a particle-mesh method, the cost of evaluating the momentum on the
grid is $\mathcal{O}(N)$, and the cost of the FFT is
$\mathcal{O}(M^2\log M)$ (where $M$ is the number of rows of
gridpoints \emph{i.e.} the total number of grid points is $M^2$),
although one could reduce this by discretising the operator used in
the norm for velocity (\emph{e.g.} the Helmholtz operator in the
examples given here) and just performing a few iterations of a method
such as Jacobi, resulting in an $\mathcal{O}(M^2)$ cost. In the case
where one is matching dense information (images or measures, for
example), then typically $N=c M^2$ (with $c>1$) and the particle-mesh
approach produces a method which is competitive with fast multipoles
methods. For the case of matching curves, typically $N = c M$ (with
$c>1$) and a good fast multipole implementation will be more
efficient. However, any mesh-based operator inversion is easily
parallelised using standard methods (as is the particle-mesh
operation), and so this could be a strength of the particle-mesh
approach in the curve case. There are also other advantages to using a
mesh, for example it is very easy to modify the operator to
investigate the behaviour with different kernels.

\section{Solution methods}
\label{solution methods}
As discussed in section \ref{discrete matching condition}, numerical
discretisation of the matching condition means that we cannot require
that the matching functional vanishes, and so we must take this into
account in the solution approach. There are two main approaches to
obtaining a numerical approximation to the optimal flow between
embedded curves:
\begin{enumerate}
\item {\bfseries Inexact matching}:
In this approach, suggested in \cite{MiTrYo03}, we 
``softly'' enforce the matching condition by adding a penalty term
to the functional \eqnref{functional}:
\begin{equation}
\label{inexact functional}
S(\MM{u}) = \int_0^1\|\MM{u}\|^2\rmd t + \frac{1}{\sigma^2}
f[\MM{v}^{\MM{Q}}].
\end{equation}
One then attempts to minimise this functional directly by applying a
descent algorithm (such as nonlinear conjugate gradients) by varying
the time series of vector fields and computing the implicit gradient
of the functional $f$.  An alternative method is to introduce the
dynamical constraint using the Lagrange multipliers $\MM{P}$ as in
section \ref{momentum}, and to solve the resulting Euler-Lagrange
equations (\ref{mom map}-\ref{P eqn}) (modified to accomodate the
penalty functional) using Newton iteration. This approach becomes
attractive when it is possible to find a good initial guess at the
solution \emph{e.g.} by specifying a path of curves between $C_A$ and
$C_B$ and approximately solving equation \eqnref{P eqn} to get an
initial guess for $\MM{P}$ given the initial guess for $\MM{Q}$. 

It is worth mentioning that the minimisation of the functional
\ref{inexact functional} becomes ill-conditioned as $\sigma\to 0$ and
so it is not always possible to match one curve onto another with the
desired accuracy (or to within the size of numerical errors, having
discretised the functional). A modification of the functional called
\emph{method of multipliers}, described in \cite{Be82}, introduces
Lagrange multiplier variables along with the penalty parameter
$\sigma$ and allows a reduction of the error in matching without
making $\sigma$ arbitrarily large. 

\item {\bfseries Minimisation by shooting}: In this approach,
  advocated in \cite{McMa2006}, we seek initial conditions for
  $\MM{P}$ such that the functional \ref{current matching} is
  minimised. The gradient of the functional with respect to the
  initial conditions for $\MM{P}$ is computed by solving the
  Euler-Lagrange equations (\ref{mom map}-\ref{P eqn}) from $t=0$ to
  $t=1$ with those initial conditions, computing the gradients of
  $f[\MM{u}^{\MM{Q}^N}]$ with respect to $\MM{Q}^N_\beta$, and
  propagating those gradients back to $\MM{P}^0_\beta$ using the
  adjoint equations (see \cite{Gu03}, for example). The problem can
  then be solved using a nonlinear gradient algorithm such as
  nonlinear conjugate gradients. As described in section \ref{matching
    condition}, $\MM{P}$ must be constrained to be normal to the curve
  in order to obtain the optimal path.

  This is the approach which we used in computing the numerical
  examples.
\end{enumerate}

\section{Numerical examples}
\label{numerics}


In this section we show a computed curve-matching calculation of two
curves in the plane. The curve was parameterised using 420 points,
with a square $128\times128$ mesh of size $2\pi\times2\pi$. The
velocity norm used was the $H^2_\alpha$-norm with $\alpha=0.4$,
discretised on the mesh using discrete Fourier transform, and the
kernel used for the current matching was the Green's function of the
$(1-\alpha^2\nabla^2)^{-2}$ operator with $\alpha=0.4$.

\begin{figure}
\begin{center}
\includegraphics[width=12cm]{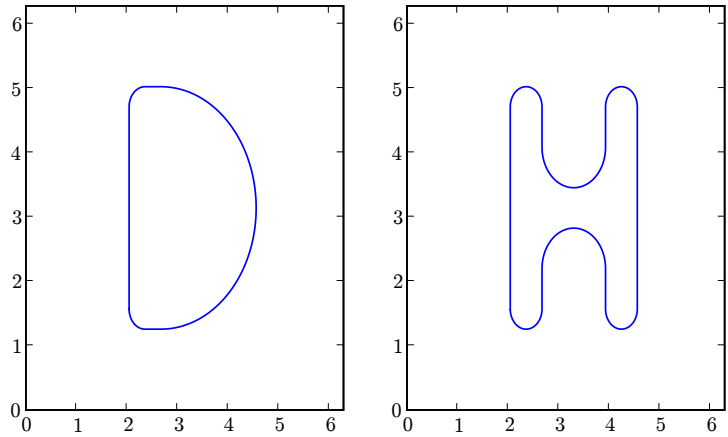}
\end{center}
\caption{\label{dandh}Test curves used for the example matching computation.}
\end{figure}

\begin{figure}
\begin{center}
\includegraphics[width=14cm]{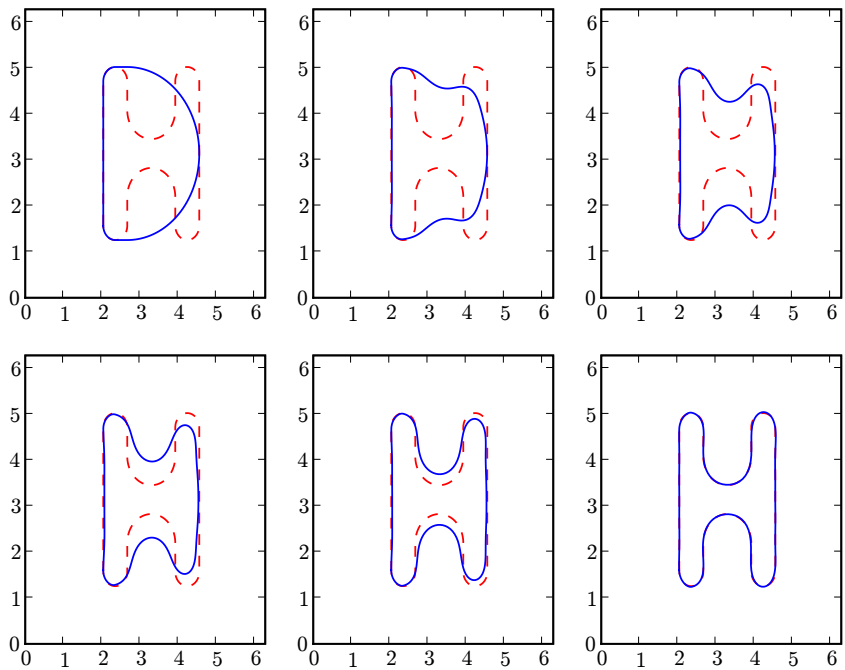}
\end{center}
\caption{\label{d2h}Snapshots of the optimal path between two test
  curves computed using the particle-mesh discretisation. The computed
  curve is plotted with a continuous line, superimposed on the target
  curve plotted with a dashed line. Top row: the curve at $t=[0,0.2,0.4]$. 
  Bottom row: the curve at $t=[0.6,0.8,1]$.  }
\end{figure}

\begin{figure}
\begin{center}
\includegraphics[width=10cm]{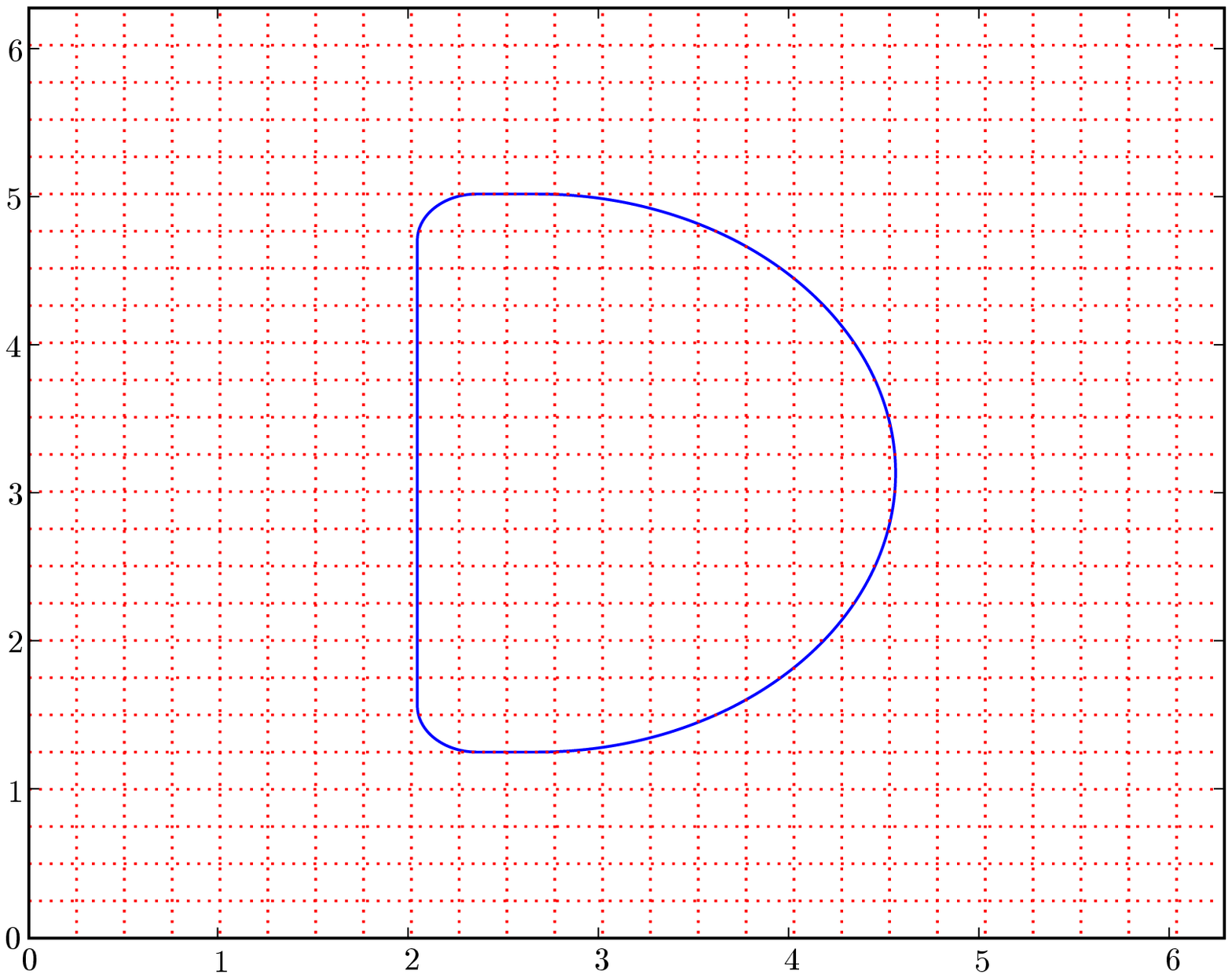}
\includegraphics[width=10cm]{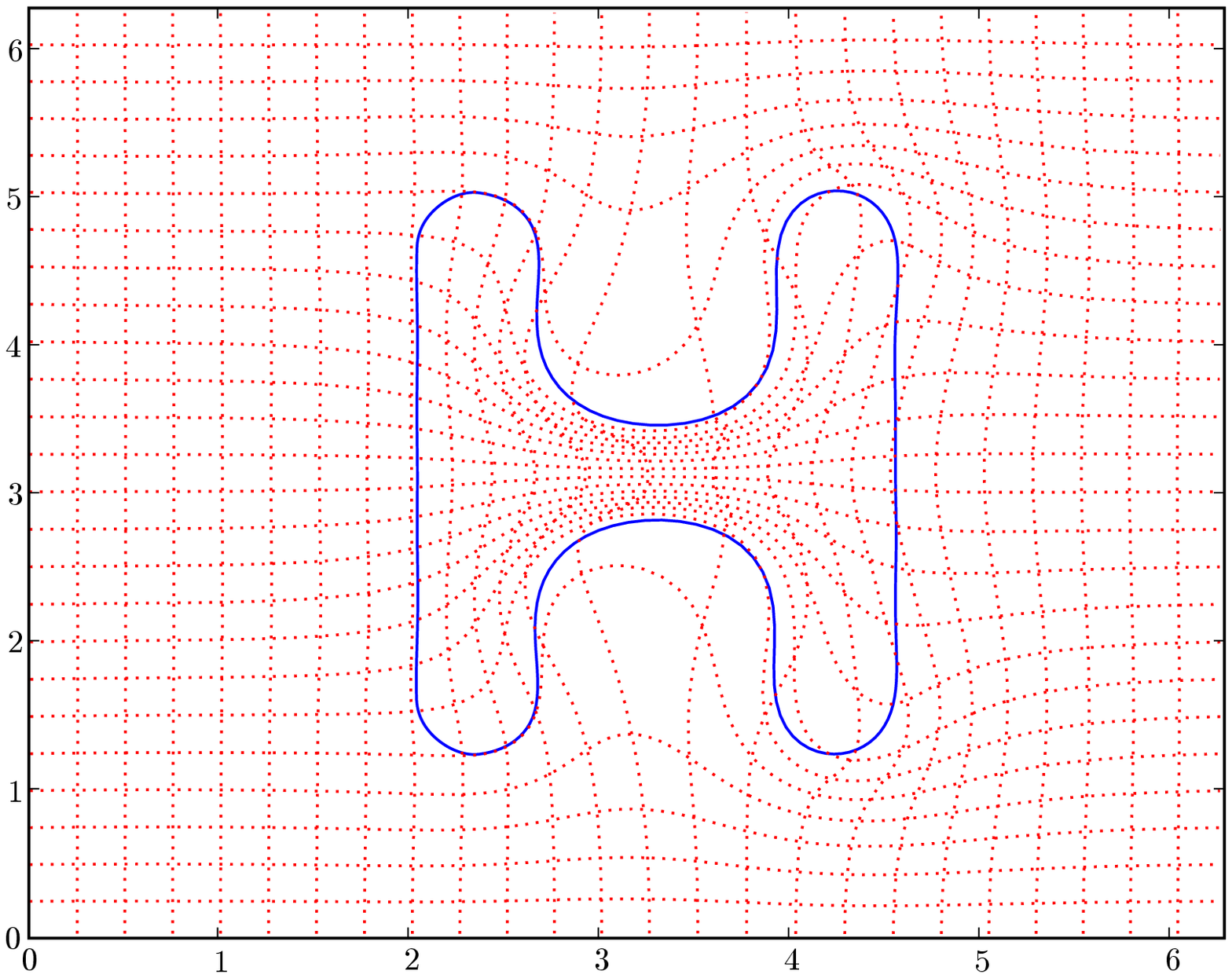}
\end{center}
\caption{\label{grid}Plots illustrating the optimal deformation map
  which maps between our two test curves. The numerical solution
  specifies a vector field which is defined everywhere, which can be
  used to transport other points which are not on the curve. This
  calculation was performed on a set of points on equispaced grid
  lines to show how space is being deformed around the curve. Top: the
  grid lines and curve before deformation. Bottom: the grid lines and
  curve after deformation. Note that the flow vector fields are far
  from divergence-free, as can be seen by inspecting the areas of the
  squares on the deformed grid. Note also that the optimal flow only
  deforms space near to the curve.}
\end{figure}

\begin{figure}\begin{center}
\includegraphics[width=12cm]{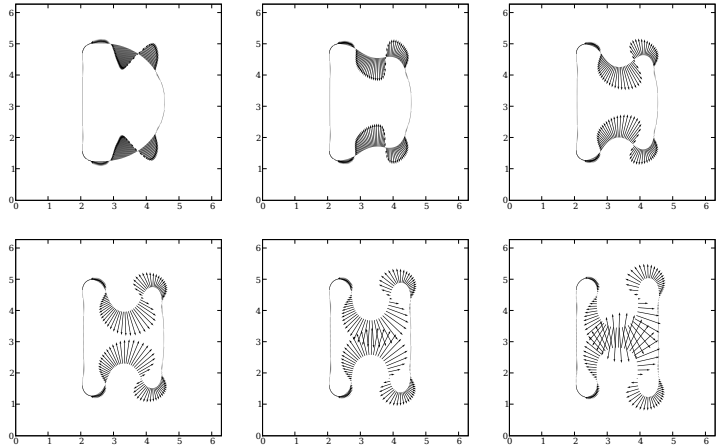}
\end{center}
\caption{\label{momentum plot}Plots showing the evolution of the momentum $\MM{P}$
  along the optimal flow. The plots are taken from the same snapshots
  as in figure \ref{d2h}, with vectors showing the direction and
  magnitude of $\MM{P}$ on the curve. Note that $\MM{P}$ remains
  normal to the curve throughout.}
\end{figure}

\begin{figure}\begin{center}
\includegraphics[width=12cm]{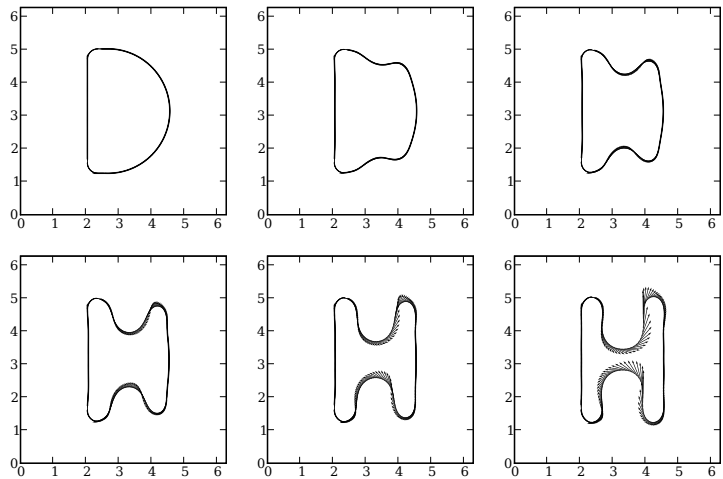}
\end{center}
\caption{\label{dQ}Plots showing $J_\beta\Delta\MM{Q}_\beta$ around
  the curve which is an approximation to $\partial\MM{Q}/\partial
  s\rmd s$ as described in section \ref{discrete relabelling}. This
  approximation remains normal to the curve throughout, with the
  increase in magnitudes showing regions of stretching in the flow.}
\end{figure}


The solution was obtained by the ``minimisation by shooting'' method
as described in section \ref{solution methods}. The minimisation was
performed using the Scientific Python \cite{scipy} {\ttfamily
  optimize.fmin\_ncg} routine which applies the Newton conjugate
gradients algorithm with the hessian computed by applying finite
differences to the gradient. The algorithm was run out until the 
functional was reduced to a value of $1.748\times 10^{-5}$ 
(with an initial value of $0.0108$).

A plot of the two curves used for the tests is given in figure
\ref{dandh}. These two curves are quite different and require large
deformations to transform one curve into the other. The calculated
path is illustrated in figure \ref{d2h} with a few snapshots of the
curve during the transformation at various times. The effect of the
deforming flow on the surrounding space is illustrated in figure
\ref{grid}.

To illustrate the results of section \ref{discrete relabelling}, a
plot showing momentum vectors at various times is given in figure
\ref{momentum plot} which suggests that the momentum stays normal to
the curve. This is confirmed by figure \ref{dQ}, showing the evolution of
$J_\beta\Delta\MM{Q}_\beta(0)$ which is an approximation to
$\partial\MM{Q}/\partial s \rmd s$.
The quantity
$\sum_\beta\MM{P}_\beta\cdot J_\beta\Delta\MM{Q}_\beta$
remains within round-off error of zero throughout, confirming the results
of section \ref{discrete relabelling}.

\section{Summary and outlook}
\label{summary}
In this paper we introduced a new particle-mesh discretisation for
diffeomorphic matching of curves, which can also be used
for matching images. In this method the vector fields used to
transport the curves are represented on a fixed mesh, whilst the
curves themselves are represented by a finite set of moving
points. Since the discrete equations arise from discretising an action
principle, they are variational integrators which have many favourable
properties. As discussed in \cite{McMa2006a}, whilst the benefits of
variational integrators have been established for long integrations
such as those for celestial mechanics or molecular dynamics, the
benefits for short time optimal control problems such as the curve
matching problem discussed in this paper are not so clear. There are a
number of drawbacks and benefits which need to be investigated with
further testing. In this paper we showed that the variational
integrators that arise from the particle-mesh discretisation have a
discrete form of the momentum conservation law which leads to the
curve momentum remaining normal to the curve throughout the computed
trajectory between shapes. We also noted that the integrators have a
modified Hamiltonian which is conserved over long times (but not
exponentially long since the compactly supported basis functions
$\psi_k(\MM{x})$ are not analytic) which can be interpreted as a
modified metric for the discrete equations. We then showed
illustrative examples obtained from the calculation of the optimal
trajectory between two closed curves in the plane.

The main computational challenge for diffeomorphic matching remains
the design of efficient algorithms to obtain optimal trajectories for
large datasets. Of the two solution methods described in section
\ref{solution methods}, the inexact matching method applied to the
particle-mesh discretisation may be best applied in parallel by using
a Newton-Krylov method to solve the Euler-Lagrange equations
(including the penalty term) since it is possible to obtain a good
initial guess for the optimal path by other methods. The
minimisation-by-shooting approach might be best applied using a
multilevel scheme where the problem is first solved with a small
number of particles, with more particles being introduced once the
lower dimensional problem has been solved to sufficient
accuracy. These approaches will need to be developed in order for the
discretisation to be applied to practical engineering applications, as
well as convergence studies and error estimates.

In future work we shall investigate the convergence properties of this
method in the limit as the number of points in the discretisation of
the curve goes to infinity (together with the number of mesh points),
and apply the method to investigate practical datasets.

Although we do not compute examples here, this method could also be
used for matching surfaces in three dimensions using a
current-matching approach. The main modifications are that equation
(\ref{u formula}) needs to be evaluated in three dimensions using a
tensor product of three B-splines (one for each Cartesian component),
and the parameterisation of the surface should become a triangulation
with a singular current being interpolated from the surface to a
three-dimensional mesh using the same basis functions. We will develop
and investigate such a method in future work. We will also investigate
efficiency of the particle-mesh versus the mesh-free fast multipoles
approach of \cite{VaGl2005} in numerical tests on real data.

Further developments will be to apply the particle mesh to the related
problem of metamorphosis \cite{TrYo2004}, and to problems where
quantities such as vector and tensor fields associated with images
also need to be matched together.

\nocite{CoHo2006,Co2005}
\bibliography{VPM}

\end{document}